\documentclass[12pt,a4paper,reqno]{amsart}
\usepackage{amsthm}
\usepackage{amsmath}
\usepackage{amsfonts}
\usepackage{amssymb}
\usepackage[left=2.8cm, right=2.8cm, top=3.0cm, bottom=3.0cm]{geometry}
\usepackage{color}
\usepackage[arrow, matrix, curve]{xy}
\usepackage[colorlinks=true,linkcolor=red,anchorcolor=blue,citecolor=green]{hyperref}

\def\be{\begin{equation}}
\def\ee{\end{equation}}
\def\bea{\begin{eqnarray}}
\def\eea{\end{eqnarray}}
\def\bes{\begin{eqnarray*}}
\def\ees{\end{eqnarray*}}

\def\nn{\nonumber}
\def\lb{\label}
\def\bs{\setminus}

\def\R{\mathbb{R}}
\def\C{\mathbb{C}}
\def\Z{\mathbb{Z}}

\def\N{\mathbb{N}}
\def\U{\mathbb{U}}

\def\Q{\mathbb{Q}}

\def\aa{{\alpha}}
\def\bb{{\beta}}
\def\ga{{\gamma}}

\def\ka{{\kappa}}
\def\th{{\theta}}

\def\om{{\omega}}

\def\ep{{\epsilon}}
\def\lm{{\lambda}}
\def\Lm{{\Lambda}}

\def\sg{{\sigma}}
\def\dm{{\diamond}}

\def\vf{{\varphi}}

\def\<{{\langle}}
\def\>{{\rangle}}

\def\P{{\mathcal{P}}}

\def\Nn{{\mathcal{N}}}

\def\rank{{\rm rank}}

\def\Sp{{\rm Sp}}

\def\mod{{\rm mod}}
\def\ol{\overline}
\def\td#1{\tilde{#1}}

\def\hb{\vrule height0.18cm width0.14cm $\,$}

\title[Closed geodesics on positively curved Finsler spheres]{On the minimal number of closed geodesics on positively curved Finsler spheres}

\author[Huagui Duan]{Huagui Duan}
\thanks{H. Duan is partially supported by National Key R\&D Program of China (Grant No. 2020YFA0713300) and NNSFC (Nos. 12271268 and 12361141812) and the Fundamental Research Funds for the Central Universities. D. Xie is partially supported by NNSFC (No. 12361141812).}
\address{(Huagui Duan) School of Mathematical Sciences and LPMC, Nankai University, Tianjin 300071, The People's Republic of China}
\email{duanhg@nankai.edu.cn.}

\author[Dong Xie]{Dong Xie}
\address{(Dong Xie) School of Mathematical Sciences, Nankai University,
Tianjin 300071, The People's Republic of China}
\email{1120200022@mail.nankai.edu.cn.}

\date{2024-05-19}

\begin{document}
\maketitle

\begin{abstract}
{\it In this paper, we proved that for every Finsler metric on $S^n$ $(n\ge 4)$ with reversibility $\lambda$ and flag curvature $K$ satisfying $( \frac{2n-3}{n-1})^2 (\frac{\lambda}{\lambda+1})^2<K\le 1$ and $ \lambda<\frac{n-1}{n-2} $, there exist at least $n$ prime closed geodesics on $(S^n,F)$, which solved a conjecture of Katok and Anosov for such positivley curved spheres when $n$ is even. Furthermore, if the number of closed geodesics on such positively curved Finsler $S^n$ is finite, then there exist at least $2\left[\frac{n}{2}\right]-1$ non-hyperbolic closed geodesics.}
\end{abstract}

\renewcommand{\theequation}{\thesection.\arabic{equation}}
\renewcommand{\thefigure}{\thesection.\arabic{figure}}

\baselineskip 18pt

\setcounter{figure}{0}
\setcounter{equation}{0}
\section{Introduction and main results}
\subsection{Introduction}

This paper is devoted to searching for the minimal number of distinct closed geodesics on positively curved Finsler spheres. This problem is closely related to a conjecture of Katok and Anosov, which claims that every Finsler metric on a sphere $S^{2k}$ or $S^{2k-1}$ carries at least $2k$ distinct closed geodesics. In fact, in 1973, Katok constructed certain irreversible Finsler metrics on spheres with flag curvature $K=1$, which have exactly $2k$ distinct closed geodesics. And based on these examples, in 1974 Anosov proposed the above conjecture.

A closed curve in a Finsler manifold is a closed geodesic if it is locally the shortest path connecting any two nearby points on this
curve. A closed geodesic $c:S^1=\R/\Z\to M$ is {\it prime}
if it is not a multiple covering (i.e., iteration) of any other
closed geodesics. Here the $m$-th iteration $c^m$ of $c$ is defined
by $c^m(t)=c(mt)$. The inverse curve $c^{-1}$ of $c$ is defined by
$c^{-1}(t)=c(1-t)$ for $t\in \R$. We call two prime closed geodesics
$c$ and $d$ {\it distinct} if there is no $\th\in (0,1)$ such that
$c(t)=d(t+\th)$ for all $t\in\R$. We shall omit the word {\it
	distinct} when we talk about more than one prime closed geodesic.
For a reversible Finsler (or Riemannian) metric, two closed geodesics
$c$ and $d$ are called { \it geometrically distinct} if $
c(S^1)\neq d(S^1)$.

For a closed geodesic $c$ on $(M,\,F)$ with $\dim M=n$, denote by $P_c$
the linearized Poincar\'{e} map of $c$. Then $P_c\in \Sp(2n-2)$ is a symplectic matrix.
We define the {\it elliptic height } $e(P_c)$
of $P_c$ to be the total algebraic multiplicity of all eigenvalues of
$P_c$ on the unit circle $\U=\{z\in\C|\; |z|=1\}$ in the complex plane
$\C$. Since $P_c$ is symplectic, $e(P_c)$ is even and $0\le e(P_c)\le 2n-2$.
Then $c$ is called {\it hyperbolic} if $e(P_c)=0$; {\it elliptic} if $e(P_c)=2n-2$;
{\it non-degenerate} if $1$ is not an eigenvalue of $P_c$. A Finsler metric $F$ on $M$ is called {\it bumpy} if all the
closed geodesics (including iterations) on $(M,F)$ are non-degenerate.

There is a famous conjecture in geometry which claims the existence of infinitely many
closed geodesics on a compact Riemannian manifold. This conjecture
has been proved for many cases, but sitll open for CROSS (compact rank one symmetric
spaces) except for $S^2$ (cf. \cite{Ban}, \cite{DLW1}, \cite{Fra}, \cite{GM}, \cite{Hin1}, \cite{Hin2} and therein).
It was surprising when Katok in \cite{Kat} found certain irreversible Finsler metrics on CROSS (including spheres) with only finitely
many closed geodesics and all of them are non-degenerate and elliptic (cf. \cite{Zil1}).

Based on the Katok's examples, in 1974, Anosov in \cite{Ano} proposed the following conjecture (also cf. \cite{Lon3} and \cite{BuM}), and we call the {\it Katok-Anosov conjecture}:
\bea \Nn(S^n,F)\ge 2\left[\frac{n+1}{2}\right] \quad \mbox{for any Finsler metric $F$ on}\ S^n, \lb{1.1}\eea
where, denote by $\Nn(M,F)$ the number of distinct closed geodesics on $(M,F)$, and $[a]=\max\{k\in\Z\,|\,k\le a\}$. In 2010, Bangert and Long in \cite{BaL} proved this conjecture for every Finsler $2$-dimensional sphere $(S^2, F)$. Since then, the index iteration theory of closed geodesics (cf. \cite{Bot} and \cite{Lon2}) has been applied to study the closed geodesic problem on Finsler manifolds.

When $n\ge 3$, the Katok-Anosov conjecture is still open in full of generality. About the multiplicity and stability problem of closed geodesics, two classes of typical conditions, including the positively curved condition and the non-degenerate (or bumpy) condition, have been studied widely.

In \cite{Rad4}, Rademacher has introduced the reversibility $\lambda=\lambda(M,F)$ of a compact Finsler manifold $(M,F)$, which is defined by
\bea \lambda=\max\{F(-X)\ |\ X\in TM,\ F(X)=1\}\ge 1.\nn\eea
Then Rademacher in \cite{Rad5} has obtained some results about the multiplicity and stability of closed geodesics. For example, let $F$ be a Finsler metric on $S^{n}$ with reversibility $\lm$ and flag curvature $K$ satisfying $\left(\frac{\lm}{1+\lm}\right)^2<K\le 1$, then there exist at least $n/2-1$ closed geodesics with length $<2n\pi$. If $\left(\frac{2n-4}{n-1}\right)^{2} \left(\frac{\lm}{1+\lm}\right)^2<K\le 1$ with $\lm<\frac{n-1}{n-3}$ and $ n\ge 4 $, then there exist at least $n-2$ closed geodesics. If $\frac{9\lm^2}{4(1+\lm)^2}<K\le 1$ with $\lm<2$, then there exists a closed geodesic of elliptic-parabolic, i.e., its linearized Poincar\'{e} map split into $2$-dimensional rotations and a part whose eigenvalues are $\pm 1$.

In 2012, Wang in \cite{Wan2} proved the conjecture (\ref{1.1}) for $(S^n,F)$ provided that $F$ is bumpy and its flag curvature $K$ satisfies $\left(\frac{\lambda}{1+\lambda}\right)^2<K\le 1$. Also in \cite{Wan2}, Wang showed that for every bumpy Finsler metric $F$ on $S^n$ satisfying $\frac{9\lm^2}{4(1+\lm)^2}<K\le 1$, there exist two prime elliptic closed geodesics if $\Nn(S^n,F)<+\infty$.
As a further generalization, Duan, Long and Wang in \cite{DLW2} obtained the optimal lower bound of the number of distinct closed geodesics on every CROSS $(M,F)$ (including spheres) if $F$ is bumpy and the flag curvature is non-negative.

Without the assumption of the non-degenerate (or bumpy) metric, the proof of conjecture (\ref{1.1}) is extremely hard. For example, under the assumption of the positively curved condition $\left(\frac{\lm}{1+\lm}\right)^2<K\le 1$, it can only be showed that $\Nn(S^n,F)\ge 3$ (cf. \cite{Dua1}) or (\cite{Dua2})) and $\Nn(S^n,F)\ge \left[ \frac{n+1}{2}\right]$ (cf. \cite{Wan4}), where $n\ge 3$.

\subsection{Main results}

Motivated by the conjecture (\ref{1.1}) and the above results, in this paper, we further consider the positively curved Finsler $n$-dimensional sphere $(S^n,F)$ without the bumpy assumption, and obtain some multiplicity results of closed geodesics on $(S^n,F)$. In particular, the following Theorem 1.1 solved the Katok-Anosov conjecture (\ref{1.1}) for certain positively curved spheres $(S^n,F)$ when $n$ is even, and this result is optimal due to the Katok's examples.

Next we assume that the reversibility $\lambda$ and flag curvature $K$ of $(S^n,F)$ satisfies the following condition:
\bea \left( \frac{2n-3}{n-1}\right)^2 \left(\frac{\lambda}{\lambda+1}\right)^2<K\le 1,\quad \mbox{where}\ \lambda<\frac{n-1}{n-2}.\lb{1.2}\eea

\medskip

{\bf Theorem 1.1.} {\it Every Finsler metric $F$ on $S^{n}$ $(n\ge 4)$ satisfying (\ref{1.2}) carries at least $n$ prime closed geodesics.}

\medskip

On the other hand, we also obtain some stability results of closed geodesics on spheres satisfying (\ref{1.2}) as follows.

\medskip

{\bf Theorem 1.2.} {\it Every Finsler metric $F$ on $S^{n}$ $(n\ge 4)$ satisfying (\ref{1.2}) carries at least $2\left[\frac{n}{2}\right]-1$ non-hyperbolic prime closed geodesics if $\Nn(S^n,F)<+\infty$. }

\medskip

{\bf Theorem 1.3.} {\it Every Finsler metric $F$ on $S^{n}$ $(n\ge 4)$ satisfying (\ref{1.2}) carries at least $n-2$ prime closed geodesics, whose linearized Poincar\'e map possesses at least an eigenvalue of the form $\exp(\pi i \mu)$ with an irrational $\mu$, if $\Nn(S^n,F)<+\infty$. }

\medskip

When the curvature condition is weaker than (\ref{1.2}), we can obtain the existence of at least $ n-1 $ prime closed geodesics, i.e. we obtained the following result.

\medskip

{\bf Theorem 1.4.} {\it On a Finsler sphere $(S^{n},F)$ $(n\ge 4)$ satisfying $\frac{9}{4}\left(\frac{\lambda}{\lambda+1}\right)^2<K\le 1$ and $ \lambda<2 $, there exist at least $ n-1 $ prime closed geodesics.}

%

\subsection{Outline of the proof}

The main ingredients of the proof of Theorem 1.1 include the Fadell-Rabinowitz index theory (see Lemma 2.5 and Lemma 2.6 below) from \cite{Rad3} and \cite{Wan3}, the Morse theory and the enhanced common index jump theorem (see Theorem 3.6 below) recently developed in \cite{DLW2}, which originally grown from the groundbreaking work \cite{LoZ}.

Let $(S^n,F)$ be a Finsler sphere satisfying (\ref{1.2}) and assume that there are only finitely many prime closed geodesics on $(S^n,F)$, which is denoted by $\{c_j\}_{j=1}^p$. The outline of the proof of Theorem 1.1 is as follows.

First, under the curvature condition (\ref{1.2}), Theorem 1.1 in \cite{Rad4} established the lower bound of the length of any prime closed geodesics on $(S^n,F)$, which in turn gives the lower bound of $ i(c_j^m) $ and the average index $ \hat{i}(c_j) $ for each closed geodesic $c_j$ (cf. Lemma 4.1 below).

Second, it follows from the enhanced common index jump theorem (cf. Theorem 3.6 below) that there exist an appropriately chosen $(p+1)$-tuple $(N,m_1,\cdots,m_p)\in\N^{p+1}$ such that the indices $i(c_j^{2m_j\pm m})$ of iterates around $2m_j$ of all $c_j$'s satisfy some precise equalities and estimates (cf. (\ref{5.6.0})-(\ref{5.9.0})). Then, together with those index estimates from Lemma 4.1, we can get much more useful and crucial estimates (cf. (\ref{5.15.0})-(\ref{5.20.0}) and Lemma 4.2) for later proofs.

Third, with these above preparations about index estimates, as the first step of the proof of Theorem 1.1, we established the existence of $n-1$ prime closed geodesics. More precisely, for each $i\in G_1\equiv[N-(n-2), N-1]\cap\N$ (and denoted by $\bar{G_1}\equiv[2N-(n-3),2N+(n-3)]\cap\N$), by the Fadell-Rabinowitz index theory (cf. Lemma 2.6), there exists some $j(i)\in\{1,\cdots,p\}$ and $m(i)\in\N$ such that the integer $2i+n-1$  (belonging to $\bar{G_1}$) lies between $i(c_{j(i)}^{m(i)})$ and $i(c_{j(i)}^{m(i)})+\nu(c_{j(i)}^{m(i)})$, and $c_{j(i)}^{m(i)}$ has non-trivial critical module $\ol{C}_{2i+n-1}(E, c_{j(i)}^{m(i)})\neq 0$. Through the analysis of index estimates and Morse theory, it can be deduced that $m(i)=2m_{j(i)}$ (cf. Claim 1 below).  Then the assumption $\Nn(S^n,F)<+\infty$ and Lemma 2.5 yield that the different integer $i\in G_1$ corresponds to the distinct closed geodesic $c_{j(i)}$. Thus it shows that $p\ge \# [N-(n-2), N-1]=n-2$. In addition, by a result of \cite{Dua1} or \cite{Dua2} (cf. Lemma 4.3 below), we can find one more closed geodesic $c_{j_0}$ satisfying  $\ol{C}_{q}(E, c_{j_0}^{2m_{j_0}})\neq 0$ only for $ q= 2N+n-1 $. Notice that $2N+n-1\notin \bar{G_1}$, thus $c_{j_0}$ is distinct from each of the $n-2$ closed geodesics found above. This yields $p=\Nn(S^n,F)\ge n-1$.

Finally, as the second step of the proof of Theorem 1.1, we carry out the proof by contradiction, i.e., assuming that there are only $ n-1 $ closed geodesics. In order to get a contradiction, here we consider another integer interval $G_2\equiv[N-(2n-3), N-(n-1)]\cap\N$ with empty intersection with $G_1$ (and denoted by $\bar{G_2}\equiv[2N-3(n-1)+2,2N-(n-1)]\cap\N$). Then, by index estimates and Morse theory, we can show that there holds $m(i)=2m_{j(i)}-1$ for every integer $i\in G_2$ and $\{c_{j(i)}\,|\, i\in G_2\}$ is the one-to-one correspondence to the $n-1$ distinct closed geodesics found in Step 1 (cf. Claim 2 and (\ref{6.6}) below). Then it can be deduced that the Morse-type number $M_{2N-(n-1)}$ only can be contributed by $c_{j(N-n+1)}^{2m_{j(N-n+1)}-1}$ (cf. Claim 3 and (\ref{6.9}) below). However, this is a contradiction $1=M_{2N-(n-1)}\ge b_{2N-(n-1)}=2$ by the Morse inequality, where note that $N$ can be chosen to be the multiple of $n-1$. Thus we obtain the existence of another distinct closed geodesic, i.e. $\Nn(S^n,F)\ge n$.

\medskip

In this paper, let $\N$, $\N_0$, $\Z$, $\Q$, $\R$, and $\C$ denote
the sets of natural integers, non-negative integers, integers,
rational numbers, real numbers, and complex numbers respectively.
We use only singular homology modules with $\Q$-coefficients.
For an $S^1$-space $X$, we denote by $\overline{X}$ the quotient space $X/S^1$.
We define the functions
\be \left\{\begin{array}{ll}
[a]=\max\{k\in\Z\,|\,k\le a\}, &
E(a)=\min\{k\in\Z\,|\,k\ge a\} , \cr
\varphi(a)=E(a)-[a], &
\{a\}=a-[a].\cr
\end{array}\right. \lb{1.3}\ee
Especially, $\varphi(a)=0$ if $ a\in\Z\,$, and $\varphi(a)=1$ if $
a\notin\Z\,$.

\setcounter{figure}{0}
\setcounter{equation}{0}

\section{Critical point theory and Fadell-Rabinowitz index theory}

\subsection{Critical point theory for closed geodesics}

On a compact Finsler manifold $(M,F)$, we choose an auxiliary Riemannian
metric. This endows the space $\Lambda=\Lambda M$ of $H^1$-maps
$\gamma:S^1\rightarrow M$ with a natural Riemannian Hilbert manifold structure
on which the group $S^1=\R/\Z$ acts continuously
by isometries, cf. \cite{Kli2}. This action is
defined by translating the parameter, i.e.,
$ (s\cdot\gamma)(t)=\gamma(t+s)$
for all $\gamma\in\Lm$ and $s,t\in S^1$.
The Finsler metric $F$ defines an energy functional $E$ and a length
functional $L$ on $\Lambda$ by
\be E(\gamma)=\frac{1}{2}\int_{S^1}F(\dot{\gamma}(t))^2dt,
\quad L(\gamma) = \int_{S^1}F(\dot{\gamma}(t))dt. \lb{2.1}\ee
Both functionals are invariant under the $S^1$-action.
The functional $E$ is $C^{1, 1}$ on $\Lm$ and
satisfies the Palais-Smale condition (cf. \cite{Mer1}). Thus we can apply the
deformation theorems in \cite{Cha} and \cite{MaW}.
The critical points
of $E$ of positive energies are precisely the closed geodesics $c:S^1\to M$
of the Finsler structure. If $c\in\Lambda$ is a closed geodesic then $c$ is
a regular curve, i.e., $\dot{c}(t)\not= 0$ for all $t\in S^1$, and this implies
that the second differential $E''(c)$ of $E$ at $c$ exists.
As usual we define the index $i(c)$ of $c$ as the maximal dimension of
subspaces of $T_c \Lambda$ on which $E^{\prime\prime}(c)$ is negative definite, and the
nullity $\nu(c)$ of $c$ so that $\nu(c)+1$ is the dimension of the null
space of $E^{\prime\prime}(c)$.

For $m\in\N$ we denote the $m$-fold iteration map
$\phi^m:\Lambda\rightarrow\Lambda$ by \be \phi^m(\ga)(t)=\ga(mt)
\qquad \forall\,\ga\in\Lm, t\in S^1. \lb{2.2}\ee We also use the
notation $\phi^m(\gamma)=\gamma^m$. For a closed geodesic $c$, the
average index is defined by \be
\hat{i}(c)=\lim_{m\rightarrow\infty}\frac{i(c^m)}{m}. \lb{2.3}\ee

If $\gamma\in\Lambda$ is not constant then the multiplicity
$m(\gamma)$ of $\gamma$ is the order of the isotropy group $\{s\in
S^1\mid s\cdot\gamma=\gamma\}$. If $m(\gamma)=1$ then $\gamma$ is
called {\it prime}. Hence $m(\gamma)=m$ if and only if there exists a
prime curve $\tilde{\gamma}\in\Lambda$ such that
$\gamma=\tilde{\gamma}^m$.

In this paper for $\ka\in \R$ and a closed geodesic $c$ we denote by
\be \Lm^{\ka}=\{d\in \Lm\,|\,E(d)\le \ka\},\qquad \Lm(c)=\{\ga\in\Lm\mid E(\ga)<E(c)\}. \lb{2.4}\ee

We call a closed geodesic satisfying the isolation condition, if
the following holds:

\medskip

{\bf (Iso)} {\it For all $m\in\N$ the orbit $S^1\cdot c^m$ is an
	isolated critical orbit of $E$.}

\medskip

Note that if the number of prime closed geodesics on a Finsler manifold
is finite, then all the closed geodesics satisfy (Iso).

Using singular homology with rational
coefficients we consider the following critical $\Q$-module of a closed geodesic
$c\in\Lambda$:
\be \overline{C}_*(E,c)
= H_*\left((\Lm(c)\cup S^1\cdot c)/S^1,\Lm(c)/S^1\right). \lb{2.5}\ee

\medskip

{\bf Proposition 2.1.} (cf. Satz 6.11 of \cite{Rad2} or Proposition
3.12 of \cite{BaL}) {\it Let $c$ be a prime closed geodesic on a Finsler manifold $(M,F)$ satisfying (Iso). Then we have
\bea
\overline{C}_q( E,c^m)
&&\equiv H_q\left((\Lm(c^m)\cup S^1\cdot c^m)/S^1, \Lm(c^m)/S^1\right)\nn\\
&&= \left(H_{i(c^m)}(U_{c^m}^-\cup\{c^m\},U_{c^m}^-)
\otimes H_{q-i(c^m)}(N_{c^m}^-\cup\{c^m\},N_{c^m}^-)\right)^{+\Z_m} \nn
\eea
	
(i) When $\nu(c^m)=0$, there holds
\bea \overline{C}_q( E,c^m) = \left\{
\begin{array}{ll}
\Q, &\quad {\it if}\;\; i(c^m)-i(c)\in 2\Z,\;{\it and}\; q=i(c^m),\; \cr
0, &\quad {\it otherwise}. \cr
\end{array}
\right.\nn
\eea
	
(ii) When $\nu(c^m)>0$, there holds
\bea
\overline{C}_q( E,c^m) = H_{q-i(c^m)}(N_{c^m}^-\cup\{c^m\}, N_{c^m}^-)^{(-1)^{i(c^m)-i(c)}\Z_m},\nn
\eea
where $N_{c^m}$ is a local characteristic manifold at $c^m$ and $N^-_{c^m}=N_{c^m}\cap \Lambda(c^m)$, $U_{c^m}$ is a local negative disk at $c^m$ and $U^-_{c^m}=U_{c^m}\cap \Lambda(c^m)$,
$H_{\ast}(X,A)^{\pm\Z_m} = \{[\xi]\in H_{\ast}(X,A)\,|\,T_{\ast}[\xi]=\pm [\xi]\}$ where $T$ is a generator of the $\Z_m$-action.}

\medskip

Denote by
\bea
k_j(c^m) \equiv \dim\, H_j(N_{c^m}^-\cup\{c^m\},N_{c^m}^- )^{(-1)^{i(c^m)-i(c)}\Z_m}. \lb{2.6}
\eea

Clearly the integers $k_j(c^m)$ equal to
$0$ when $j<0$ or $j>\nu(c^m)$ and can take only values $0$ or $1$
when $j=0$ or $j=\nu(c^m)$.

\medskip

{\bf Proposition 2.2.} (cf. Satz 6.13 of \cite{Rad2}) {\it Let $c$
be a prime closed geodesic on a Finsler manifold $(M, F)$ satisfying (Iso). For any $m\in\N$, we have
	
(i) If $k_0(c^m)=1$, there holds $k_j(c^m)=0$ for $1\le j\le \nu(c^m)$.
	
(ii) If $k_{\nu(c^m)}(c^m)=1$, there holds $k_j(c^m)=0$ for $0\le j\le \nu(c^m)-1$.
	
(iii) If $k_j(c^m)\ge 1$ for some $1\le j\le \nu(c^m)-1$, there holds $k_{\nu(c^m)}(c^m)=0=k_0(c^m)$.}

\medskip

Set $\ol{\Lm}^0=\ol{\Lambda}^0M =\{{\rm constant\;point\;curves\;in\;}M\}\cong M$. Let $(X,Y)$ be a
space pair such that the Betti numbers $b_i=b_i(X,Y)=\dim
H_i(X,Y;\Q)$ are finite for all $i\in \Z$. As usual the {\it
	Poincar\'e series} of $(X,Y)$ is defined by the formal power series
$P(X, Y)=\sum_{i=0}^{\infty}b_it^i$. We need the following results on Betti numbers.

\medskip

{\bf Theorem 2.3.} (cf. \cite{Hin0}, Theorem 2.4 and Remark 2.5 of \cite{Rad1})
{\it We have the Poincar\'e series of $(\ol{\Lm}S^n,\ol{\Lm}^0S^n)$ as the following:
	
(i) When $n=2k+1$ is odd
\bea
P(\ol{\Lm}S^n,\ol{\Lm}^0S^n)(t)
=t^{n-1}\left(\frac{1}{1-t^2}+\frac{t^{n-1}}{1-t^{n-1}}\right)= t^{2k}\left(\frac{1}{1-t^2}+\frac{t^{2k}}{1- t^{2k}}\right).
\lb{2.14}
\eea
Thus for $q\in\Z$ and $l\in\N_0$, we have
\bea
{b}_q &=& {b}_q(\ol{\Lm}S^n,\ol{\Lm}^0 S^n)=\rank H_q(\ol{\Lm} S^n,\ol{\Lm}^0 S^n )\nn\\
	&=& \;\;\left\{
\begin{array}{ll}
	2,&\quad {\it if}\quad q\in \{4k+2l,\quad l=0\;\mod\; k\}, \cr
1,&\quad {\it if}\quad q\in \{2k\}\cup\{2k+2l,\quad l\neq 0\;\mod\; k\}, \cr
0 &\quad {\it otherwise}. \cr
\end{array}
\right. \lb{2.15}
\eea
	
(ii) When $n=2k$ is even
\bea
P(\ol{\Lm}S^n,\ol{\Lm}^0S^n)(t)
&=&t^{n-1}\left(\frac{1}{1-t^2}+\frac{t^{n(m+1)-2}}{1-t^{n(m+1)-2}}\right) \frac{1-t^{nm}}{1-t^n}\nn\\
&=& t^{2k-1}\left(\frac{1}{1-t^2}+\frac{t^{4k-2}}{1- t^{4k-2}}\right),\lb{2.16}
\eea
where $m=1$ by Theorem 2.4 of \cite{Rad1}. Thus for $q\in\Z$ and $l\in\N_0$, we have}
\bea
{b}_q &=& {b}_q(\ol{\Lm}S^n,\ol{\Lm}^0 S^n)=\rank H_q(\ol{\Lm} S^n,\ol{\Lm}^0 S^n )\nn\\
&=& \;\;\left\{
\begin{array}{ll}
	2,&\quad {\it if}\quad q\in \{6k-3+2l,\quad l=0\;\mod\; 2k-1\}, \cr
1,&\quad {\it if}\quad q\in \{2k-1\}\cup\{2k-1+2l,\quad l\neq 0\;\mod\; 2k-1\}, \cr
0 &\quad {\it otherwise}. \cr
\end{array}
\right. \lb{2.17}
\eea

We have the following version of the Morse inequality.

\medskip

{\bf Theorem 2.4.} (cf. Theorem 6.1 of \cite{Rad2}) {\it Suppose that there exist
	only finitely many prime closed geodesics $\{c_j\}_{1\le j\le p}$ on $(M, F)$,
	and $0\le a<b\le \infty$.
	Define for each $q\in\Z$,
	\bea
	{M}_q(\ol{\Lm}^b,\ol{\Lm}^a)
	&=& \sum_{1\le j\le p,\;a<E(c^m_j)<b}\rank{\ol{C}}_q(E, c^m_j ) \nn\\
	{b}_q(\ol{\Lm}^{b},\ol{\Lm}^{a})
	&=& \rank H_q(\ol{\Lm}^{b},\ol{\Lm}^{a}). \nn\eea
	Then there holds }
\bea
M_q(\ol{\Lm}^{b},\ol{\Lm}^{a}) &-& M_{q-1}(\ol{\Lm}^{b},\ol{\Lm}^{a})
+ \cdots +(-1)^{q}M_0(\ol{\Lm}^{b},\ol{\Lm}^{a}) \nn\\
&\ge& b_q(\ol{\Lm}^{b},\ol{\Lm}^{a}) - b_{q-1}(\ol{\Lm}^{b},\ol{\Lm}^{a})
+ \cdots + (-1)^{q}b_0(\ol{\Lm}^{b},\ol{\Lm}^{a}), \lb{2.18}\\
{M}_q(\ol{\Lm}^{b},\ol{\Lm}^{a}) &\ge& {b}_q(\ol{\Lm}^{b},\ol{\Lm}^{a}).\lb{2.19}
\eea

\subsection{Fadell-Rabinowitz index theory}

Next we recall the Fadell-Rabinowitz index
in a relative version due to \cite{Rad3}.
Let $X$ be an $S^1$-space, $A\subset X$ a closed $S^1$-invariant subset.
Note that the cup product defines a homomorphism
\be H^\ast_{S^1}(X)\otimes H^\ast_{S^1}(X,\; A)\rightarrow H^\ast_{S^1}(X,\;A):
\quad (\zeta,\;z)\rightarrow \zeta\cup z,\lb{2.20}\ee
where $H_{S^1}^\ast$ is the $S^1$-equivariant cohomology with
rational coefficients in the sense of Borel
(cf. Chapter IV of \cite{Bor}).
Fix a characteristic class $\eta\in H^2(CP^\infty)$. Let
$f^\ast: H^\ast(CP^\infty)\rightarrow H_{S^1}^\ast(X)$ be the
homomorphism induced by a classifying map $f: X_{S^1}\rightarrow CP^\infty$.
Now for $\gamma\in H^\ast(CP^\infty)$ and $z\in H^\ast_{S^1}(X,\;A)$, let
$\gamma\cdot z=f^\ast(\gamma)\cup z$. Then
the order $\mbox{ord}_\eta(z)$ with respect to $\eta$ is defined by
\be \mbox{ord}_\eta(z)=\inf\{k\in\N\cup\{\infty\}\;|\;\eta^k\cdot z= 0\}.\lb{2.21}\ee
By Proposition 3.1 of \cite{Rad3}, there is an element
$z\in H^{n+1}_{S^1}(\Lambda,\;\Lambda^0)$ of infinite order, i.e., $\mbox{ord}_\eta(z)=\infty$.
For $\kappa\ge0$, we denote by
$j_\kappa: (\Lambda^\kappa, \;\Lambda^0)\rightarrow (\Lambda \;\Lambda^0)$
the natural inclusion and define the function $d_z:\R^{\ge 0}\rightarrow \N\cup\{\infty\}$:
\be d_z(\kappa)=\mbox{ord}_\eta(j_\kappa^\ast(z)).\lb{2.22}\ee
Denote by $d_z(\kappa-)=\lim_{\epsilon\searrow 0}d_z(\kappa-\epsilon)$,
where $t\searrow a$ means $t>a$ and $t\to a$.

By the Section 5 of \cite{Rad3}, the function $d_z$ is non-decreasing and $\lim_{\lambda\searrow \kappa}d_z(\lambda)=d_z(\kappa).$ Each discontinuous point of $d_z$ is a critical value of the energy functional $E$. If $d_z(\kappa)-d_z(\kappa-)\ge 2$, then there are infinitely many prime closed geodesics $c$ with energy $\kappa$. So, for each $i\ge 1$, we define $$\kappa_i=\inf\{\delta\in\R\;|\: d_z(\delta)\ge i\}.$$ Then we have the following:

\medskip

{\bf Lemma 2.5.} (cf. Lemma 2.3 of \cite{Wan3}) {\it Suppose there are only finitely many prime closed
	geodesics on $(S^n,\, F)$. Then each $\kappa_i$ is a critical value of $E$.
	If $\kappa_i=\kappa_j$ for some $i<j$, then there are infinitely many prime
	closed geodesics on $(S^n,\, F)$. }

\medskip

{\bf Lemma 2.6.} (cf. Lemma 2.4 of \cite{Wan3}) {\it Suppose there are only finitely many prime closed
	geodesics on $(S^n,\, F)$. Then for every $i\in\N$, there exists a
	closed geodesic $c$ on $(S^n,\, F)$ such that
	\bea E(c)=\kappa_i,\quad
	\ol{C}_{2i+\dim(z)-2}(E, c)\neq 0.\lb{2.24}\eea}

{\bf Definition 2.7.} {\it A prime closed geodesic $c$ is
	$(m, i)$-variationally visible, if there exist some $m, i\in\N$
	such that (\ref{2.24}) holds for $c^m$ and $\kappa_i$. We call $c$
	infinitely variationally visible, if there exist
	infinitely many $m, i\in\N$ such that $c$ is $(m, i)$-variationally visible.
	Denote by $\mathcal{V}_\infty(S^n, F)$ the set of infinitely
	variationally visible closed geodesics.}

\medskip

{\bf Theorem 2.8.} (cf. Theorem 2.6 of \cite{Wan3}) {\it Suppose there are only finitely many prime closed
	geodesics on $(S^n,\, F)$. Then for any $c\in\mathcal{V}_\infty(S^n, F)$,
	there holds
	\be \frac{\hat i(c)}{L(c)}=2\sigma.\lb{2.25}
	\ee
	where $\sigma=\lim\inf_{i\rightarrow\infty}i/\sqrt{2\kappa_i}=
	\lim\sup_{i\rightarrow\infty}{i}/{\sqrt{2\kappa_i}}$.
}

\setcounter{equation}{0}
\section{Index iteration theory for closed geodesics}

In this section, we recall briefly the index theory for symplectic paths
developed by Long and his coworkers, and all the details can be found in \cite{Lon2}.

Let $c$ be a closed geodesic on an orientable Finsler manifold $M=(M,\,F)$.
Denote the linearized Poincar\'e map of $c$ by $P_c\in\Sp(2n-2)$.
Then $P_c$ is a symplectic matrix.
Note that the index iteration formulae in \cite{Lon1} of 2000 (cf. Chap. 8 of
\cite{Lon2}) work for Morse indices of iterated closed geodesics (cf. \cite{LLo}, Chap. 12 of \cite{Lon2}). Since every closed geodesic
on $M$ is orientable. Then by Theorem 1.1 of \cite{Liu1}
of Liu (cf. also \cite{Wil}), the initial Morse index of a closed geodesic
$c$ on $M$ coincides with the index of a
corresponding symplectic path introduced by Conley, Zehnder and Long
in 1984-1990 (cf. \cite{Lon2}).

As usual, the symplectic group $\Sp(2n)$ is defined by
$$ \Sp(2n) = \{M\in {\rm GL}(2n,\R)\,|\,M^TJM=J\}, $$
whose topology is induced from that of $\R^{4n^2}$. For $\tau>0$ we are interested
in $\Sp(2n)$:
$$ \P_{\tau}(2n) = \{\ga\in C([0,\tau],\Sp(2n))\,|\,\ga(0)=I_{2n}\}, $$
which is equipped with the topology induced from that of $\Sp(2n)$. The
following real function was introduced in \cite{Lon1}:
$$ D_{\om}(M) = (-1)^{n-1}\ol{\om}^n\det(M-\om I_{2n}), \qquad
\forall \om\in\U,\, M\in\Sp(2n). $$
Thus for any $\om\in\U$ the following hypersurface in $\Sp(2n)$ is
defined in \cite{Lon1}:
$$ \Sp(2n)_{\om}^0 = \{M\in\Sp(2n)\,|\, D_{\om}(M)=0\}. $$
For any $M\in \Sp(2n)_{\om}^0$, we define a co-orientation of $\Sp(2n)_{\om}^0$
at $M$ by the positive direction $\frac{d}{dt}Me^{t\ep J}|_{t=0}$ of
the path $Me^{t\ep J}$ with $0\le t\le 1$ and small $\ep>0$. Let
\bea
\Sp(2n)_{\om}^{\ast} &=& \Sp(2n)\bs \Sp(2n)_{\om}^0, \nn\\
\P_{\tau,\om}^{\ast}(2n) &=&
\{\ga\in\P_{\tau}(2n)\,|\,\ga(\tau)\in\Sp(2n)_{\om}^{\ast}\}, \nn\\
\P_{\tau,\om}^0(2n) &=& \P_{\tau}(2n)\bs \P_{\tau,\om}^{\ast}(2n). \nn\eea
For any two continuous arcs $\xi$ and $\eta:[0,\tau]\to\Sp(2n)$ with
$\xi(\tau)=\eta(0)$, we define
$$ \eta\ast\xi(t) = \left\{
\begin{array}{ll}
\xi(2t), & \quad {\rm if}\;0\le t\le \tau/2, \cr
\eta(2t-\tau), & \quad {\rm if}\; \tau/2\le t\le \tau. \cr
\end{array}
\right. $$
Given any two $2m_k\times 2m_k$ matrices of square block form
$M_k=\left(
\begin{array}{cc}
	A_k&B_k\cr
C_k&D_k\cr
\end{array}
\right)$ with $k=1, 2$,
as in \cite{Lon2}, the $\;\diamond$-product of $M_1$ and $M_2$ is defined by
the following $2(m_1+m_2)\times 2(m_1+m_2)$ matrix $M_1\diamond M_2$:
$$ M_1\diamond M_2=\left(
\begin{array}{cccc}
A_1& 0&B_1& 0\cr
0&A_2& 0&B_2\cr
C_1& 0&D_1& 0\cr
0&C_2& 0&D_2\cr
\end{array}
\right). \nn$$ 
Denote by $M^{\diamond k}$ the $k$-fold $\diamond$-product $M\diamond\cdots\diamond M$. Note
that the $\diamond$-product of any two symplectic matrices is symplectic. For any two
paths $\ga_j\in\P_{\tau}(2n_j)$ with $j=0$ and $1$, let
$\ga_0\diamond\ga_1(t)= \ga_0(t)\diamond\ga_1(t)$ for all $t\in [0,\tau]$.

A special path $\xi_n$ is defined by
\be \xi_n(t) = \left(
\begin{array}{cc}
2-\frac{t}{\tau} & 0 \cr
0 & (2-\frac{t}{\tau})^{-1}\cr
\end{array}
\right)^{\diamond n}
\qquad {\rm for}\;0\le t\le \tau. \lb{3.1}\ee

{\bf Definition 3.1.} (cf. \cite{Lon1}, \cite{Lon2}) {\it For any $\om\in\U$ and $M\in \Sp(2n)$, define
\be
\nu_{\om}(M)=\dim_{\C}\ker_{\C}(M - \om I_{2n}). \lb{3.2}
\ee
For any $\tau>0$ and $\ga\in \P_{\tau}(2n)$, define
\be
\nu_{\om}(\ga)= \nu_{\om}(\ga(\tau)). \lb{3.3}
\ee
	
If $\ga\in\P_{\tau,\om}^{\ast}(2n)$, define
\be
i_{\om}(\ga) = [\Sp(2n)_{\om}^0: \ga\ast\xi_n], \lb{3.4}
\ee
where the right hand side of (\ref{3.4}) is the usual homotopy intersection number, and the orientation of $\ga\ast\xi_n$ is its positive time direction under homotopy with fixed end points.
	
If $\ga\in\P_{\tau,\om}^0(2n)$, we let $\mathcal{F}(\ga)$ be the set of all open neighborhoods of $\ga$ in $\P_{\tau}(2n)$, and define
\be
i_{\om}(\ga) = \sup_{U\in\mathcal{F}(\ga)}\inf\{i_{\om}(\beta)\,|\,
\beta\in U\cap\P_{\tau,\om}^{\ast}(2n)\}.
\lb{3.5}
\ee
Then
$$
(i_{\om}(\ga), \nu_{\om}(\ga)) \in \Z\times \{0,1,\ldots,2n\}, $$
is called the index function of $\ga$ at $\om$. }

\medskip

For any symplectic path $\ga\in\P_{\tau}(2n)$ and $m\in\N$, we
define its $m$-th iteration $\ga^m:[0,m\tau]\to\Sp(2n)$ by
\be \ga^m(t) = \ga(t-j\tau)\ga(\tau)^j, \qquad
{\rm for}\quad j\tau\leq t\leq (j+1)\tau,\;j=0,1,\ldots,m-1.
\lb{3.6}\ee
We still denote the extended path on $[0,+\infty)$ by $\ga$.

\medskip

{\bf Definition 3.2.} (cf. \cite{Lon1}, \cite{Lon2}) {\it For any $\ga\in\P_{\tau}(2n)$,
we define
\be
(i(\ga,m), \nu(\ga,m)) = (i_1(\ga^m), \nu_1(\ga^m)), \qquad \forall m\in\N.
\lb{3.7}
\ee
The average index $\hat{i}(\ga,m)$ per $m\tau$ for $m\in\N$ is defined by
\be
\hat{i}(\ga,m) = \lim_{k\to +\infty}\frac{i(\ga,mk)}{k}. \lb{3.8}
\ee
For $M\in\Sp(2n)$ and $\om\in\U$, the {\it splitting numbers} $S_M^{\pm}(\om)$ of $M$ at $\om$ are defined by
\be
S_M^{\pm}(\om) = \lim_{\ep\to 0^+}i_{\om\exp(\pm\sqrt{-1}\ep)}(\ga) - i_{\om}(\ga),\lb{3.9}
\ee
for any path $\ga\in\P_{\tau}(2n)$ satisfying $\ga(\tau)=M$.}

\medskip

For a given path $\gamma\in \P_{\tau}(2n)$, we deform it to a new path $\eta$ in $\P_{\tau}(2n)$ so that
\begin{equation}
i_1(\gamma^m)=i_1(\eta^m),\quad \nu_1(\gamma^m)=\nu_1(\eta^m), \quad
\forall m\in {\N}, \label{3.10}
\end{equation}
and that $(i_1(\eta^m),\nu_1(\eta^m))$ is easy enough to compute. This
leads to finding homotopies $\delta:[0,1]\times[0,\tau]\to {\rm Sp}(2n)$
starting from $\gamma$ in $\P_{\tau}(2n)$ and keeping the end
points of the homotopy always stay in a certain suitably chosen maximal
subset of ${\rm Sp}(2n)$ so that (\ref{3.10}) always holds. In fact, this
set was first discovered in \cite{Lon1} as the path connected component
$\Omega^0(M)$ containing $M=\gamma(\tau)$ of the set
\begin{eqnarray}
\Omega(M)=\{N\in{\rm Sp}(2n)\,&|&\,\sigma(N)\cap{\U}=\sigma(M)\cap{\U}\;
{\rm and}\; \nonumber\\
&&\qquad \nu_{\lambda}(N)=\nu_{\lambda}(M),\;\forall\,
\lambda\in\sigma(M)\cap{\U}\}. \label{3.11}
\end{eqnarray}
Here $\Omega^0(M)$ is called the {\it homotopy component} of $M$ in
${\rm Sp}(2n)$.

In \cite{Lon1} and \cite{Lon2}, the following symplectic matrices were introduced
as {\it basic normal forms}:
\begin{eqnarray}
D(\lambda)=\left(
\begin{array}{cc}
\lm & 0\cr
0 & \lm^{-1}\cr
\end{array}
\right), &\quad& \lm=\pm 2,\lb{3.12}\\
N_1(\lm,b) = \left(
\begin{array}{cc}
\lm & b\cr
0 & \lm\cr
\end{array}
\right), &\quad& \lm=\pm 1, b=\pm1, 0, \lb{3.13}\\
R(\th)=\left(
\begin{array}{cc}
\cos\th & -\sin\th\cr
\sin\th & \cos\th\cr
\end{array}
\right), &\quad& \th\in (0,\pi)\cup(\pi,2\pi),\lb{3.14}\\
N_2(\om,b)= \left(
\begin{array}{cc}
R(\th) & b\cr
0 & R(\th)\cr
\end{array}
\right), &\quad& \th\in (0,\pi)\cup(\pi,2\pi),
\lb{3.15}
\end{eqnarray}
where $b=\left(
\begin{array}{cc}
b_1 & b_2\cr
b_3 & b_4\cr
\end{array}
\right)$ with $b_i\in\R$ and $b_2\not=b_3$.
We call $ N_2(\omega, b)$
{\it trivial} if $(b_2-b_3)\sin\theta>0$, $ N_2(\omega, b)$
{\it non-trivial} if $(b_2-b_3)\sin\theta<0$.

Splitting numbers possess the following properties:

\medskip

{\bf Lemma 3.3.} (cf. \cite{Lon1}, Lemma 9.1.5 and List 9.1.12 of \cite{Lon2})
{\it For $M\in\Sp(2n)$ and $\om\in\U$, there hold
\bea
S_M^{\pm}(\om) &=& 0, \qquad {\it if}\;\;\om\not\in\sg(M). \lb{3.16}\\
S_{N_1(1,a)}^+(1) &=& \left\{
\begin{array}{ll}
1, &\quad {\rm if}\;\; a\ge 0, \cr
0, &\quad {\rm if}\;\; a< 0. \cr
\end{array}
\right. \lb{3.17}
\eea
	
	For any $M_i\in\Sp(2n_i)$ with $i=0$ and $1$, there holds }
\be S^{\pm}_{M_0\diamond M_1}(\om) = S^{\pm}_{M_0}(\om) + S^{\pm}_{M_1}(\om),
\qquad \forall\;\om\in\U. \lb{3.18}\ee

The following is the precise index iteration formulae for symplectic paths, which is due to Long (cf. Theorem 8.3.1 and Corollary 8.3.2 of \cite{Lon2}).

\medskip

{\bf Theorem 3.4.} {\it Let $\gamma\in\P_\tau(2n)$, then there exists a path $f\in C([0,1],\Omega^0(\gamma(\tau))$ such that $f(0)=\gamma(\tau)$ and
\bea
f(1)=&&N_1(1,1)^{\diamond p_-} \diamond I_{2p_0}\diamond
N_1(1,-1)^{\diamond p_+}
\diamond N_1(-1,1)^{\diamond q_-} \diamond (-I_{2q_0})\diamond
N_1(-1,-1)^{\diamond q_+}\nn\\
&&\diamond R(\theta_1)\diamond\cdots\diamond R(\theta_r)
\diamond N_2(\omega_1, u_1)\diamond\cdots\diamond N_2(\omega_{r_*}, u_{r_*}) \nn\\
&&\diamond N_2(\lm_1, v_1)\diamond\cdots\diamond N_2(\lm_{r_0}, v_{r_0}) \diamond M_0 \lb{3.19}
\eea
where $ N_2(\omega_j, u_j) $s are non-trivial and $ N_2(\lm_j, v_j)$s are trivial basic normal forms; $\sigma (M_0)\cap U=\emptyset$; $p_-$, $p_0$, $p_+$, $q_-$, $q_0$, $q_+$, $r$, $r_*$ and $r_0$ are non-negative integers; $\omega_j=e^{\sqrt{-1}\alpha_j}$, $\lambda_j=e^{\sqrt{-1}\beta_j}$; $\theta_j$, $\alpha_j$, $\beta_j$
$\in (0, \pi)\cup (\pi, 2\pi)$; these integers and real numbers
are uniquely determined by $\gamma(\tau)$. Then using the
functions defined in (\ref{1.3})
\bea
i(\gamma, m)&=&m(i(\gamma, 1)+p_-+p_0-r)+2\sum_{j=1}^r E\left(\frac{m\theta_j}{2\pi}\right)-r -p_--p_0\nn\\&&-\frac{1+(-1)^m}{2}(q_0+q_+)+2\left(\sum_{j=1}^{r_*}\varphi\left(\frac{m\alpha_j}{2\pi}\right)-r_*\right).
\lb{3.20}\\
\nu(\gamma, m)&=&\nu(\gamma, 1)+\frac{1+(-1)^m}{2}(q_-+2q_0+q_+)+2(r+r_*+r_0)\nn\\
&&-2\left(\sum_{j=1}^{r}\varphi\left(\frac{m\theta_j}{2\pi}\right)+
\sum_{j=1}^{r_*}\varphi\left(\frac{m\alpha_j}{2\pi}\right)
+\sum_{j=1}^{r_0}\varphi\left(\frac{m\beta_j}{2\pi}\right)\right)\lb{3.21}\\
\hat i(\gamma, 1)&=&i(\gamma, 1)+p_-+p_0-r+\sum_{j=1}^r \frac{\theta_j}{\pi}.\lb{3.22}
\eea

In addition, $i(\gamma, 1)$ is odd if $f(1)=N_1(1, 1)$, $I_2$,
$N_1(-1, 1)$, $-I_2$, $N_1(-1, -1)$ and $R(\theta)$; $i(\gamma, 1)$ is even if $f(1)=N_1(1, -1)$ and $ N_2(\omega, b)$; $i(\gamma, 1)$ can be any integer if $\sigma (f(1)) \cap \U=\emptyset$.}

\medskip

We have the following properties in the index iteration theory.

\medskip

{\bf Theorem 3.5. } (cf. Theorem 2.2 of \cite{LoZ}) {\it Let $\gamma\in\P_\tau(2n)$ and $M = \gamma(\tau)$.
Then for any $m\in\N$, there holds
$$
\nu(\gamma, m)-\frac{e(M)}{2}\le i(\gamma, m+1)-i(\gamma, m)-i(\gamma, 1) \le \nu(\gamma, 1)-\nu(\gamma, m+1)+\frac{e(M)}{2},
$$
where $e(M)$ is the elliptic hight defined in Section 1.}

\medskip

The common index jump theorem (cf. Theorem 4.3 of \cite{LoZ}) for symplectic paths has become
one of the main tools in studying periodic orbits in Hamiltonian and symplectic dynamics. In 2016, the following enhanced common index jump theorem has been
obtained by Duan, Long and Wang in \cite{DLW2}.

\medskip

{\bf Theorem 3.6.} (cf. Theorem 3.5 of \cite{DLW2}) {\it Let
$\gamma_k\in\mathcal{P}_{\tau_k}(2n)$ for $k=1,\cdots,q$ be a finite collection of symplectic paths. Let $M_k=\ga_k(\tau_k)$. We extend $\ga_k$ to $[0,+\infty)$ by (\ref{3.6}) inductively. Suppose
\be
\hat{i}(\ga_k,1) > 0, \qquad \forall\ k=1,\cdots,q. \lb{3.23}
\ee
	Then for any fixed integer $\bar{m}\in \N$, there exist infinitely many $(q+1)$-tuples
	$(N, m_1,\cdots,\\m_q) \in \N^{q+1}$ such that for all $1\le k\le q$ and $1\le m\le \bar{m}$, there holds
	\bea
	\nu(\ga_k,2m_k-m) &=& \nu(\ga_k,2m_k+m) = \nu(\ga_k, m),  \lb{3.24}\\
	i(\ga_k,2m_k+m) &=& 2N+i(\ga_k,m),             \lb{3.25}\\
	i(\ga_k,2m_k-m) &=& 2N-i(\ga_k,m)-2(S^+_{M_k}(1)+Q_k(m)), \lb{3.26}\\
	i(\ga_k, 2m_k)&=& 2N -(S^+_{M_k}(1)+C(M_k)-2\Delta_k),   \lb{3.27}\eea
	where $S_{M_k}^\pm(\om)$ is the splitting number of $M_k$ at $\om$ in Definition 3.2 and
	\bea
	&&C(M_k)=\sum\limits_{0<\theta<2\pi}S^-_{M_k}(e^{\sqrt{-1}\theta}),\ \Delta_k = \sum_{0<\{m_k\th/\pi\}<\delta}S^-_{M_k}(e^{\sqrt{-1}\th}),\nn\\
	&&Q_k(m) = \sum_{\theta\in(0,2\pi), e^{\sqrt{-1}\th}\in\sg(M_k),\atop
		\{\frac{m_k\th}{\pi}\}= \{\frac{m\th}{2\pi}\}=0}
	S^-_{M_k}(e^{\sqrt{-1}\th}). \lb{3.28}
	\eea
	More precisely, by (4.10), (4.40) and (4.41) in \cite{LoZ}, we have
	\bea m_k=\left(\left[\frac{N}{\bar{M}\hat i(\gamma_k, 1)}\right]+\chi_k\right)\bar{M},\quad 1\le k\le q,\lb{3.29}\eea
	where $\chi_k=0$ or $1$ for $1\le k\le q$ and $\frac{\bar{M}\theta}{\pi}\in\Z$
	whenever $e^{\sqrt{-1}\theta}\in\sigma(M_k)$ and $\frac{\theta}{\pi}\in\Q$
	for some $1\le k\le q$. Furthermore, for any fixed $M_0\in\N$, we may
	further require $M_0|N$, and for any $\epsilon>0$, we can choose $N$ and $\{\chi_k\}_{1\le k\le q}$ such that
	\bea \left|\left\{\frac{N}{\bar{M}\hat i(\gamma_k, 1)}\right\}-\chi_k\right|<\epsilon,\quad 1\le k\le q.\lb{3.30}\eea}

\setcounter{figure}{0}
\setcounter{equation}{0}
\section{Proof of main results}

\subsection{Proof of Theorem 1.1}

In this subsection, we make the following assumption and give the proof of Theorem 1.1.

\medskip

{\bf (F1)} {\it For any Finsler sphere $(S^{n},F)$ $(n\ge 4)$ satisfying $\left(\frac{2n-3}{n-1}\right)^2\left(\frac{\lm}{1 + \lm}\right)^2 < K \le 1$ and $\lm<\frac{n-1}{n-2}$, there holds $\Nn(S^n,F)<+\infty$, denoted by $\{c_j\}_{j=1}^p$.}

\medskip

For any $1\le j\le p$, we rewrite (\ref{3.19}) as follows
\bea f_j(1)
&=& N_1(1,1)^{\dm p_{j,-}}\,\dm\,I_{2p_{j,0}}\,\dm\,N_1(1,-1)^{\dm p_{j,+}}\nn\\
&&\dm\,N_1(-1,1)^{\dm q_{j,-}}\,\dm\,(-I_{2q_{j,0}})\,\dm\,N_1(-1,-1)^{\dm q_{j,+}} \nn\\
&& \dm\,R(\th_{j,1})\,\dm\,\cdots\,\dm\,R(\th_{j,r_{j,1}})\,\dm\,R(\td{\th}_{j,1})\,\dm\,\cdots\,\dm\,R(\td{\th}_{j,r_{j,2}})\nn\\
&& \dm\,N_2(e^{\sqrt{-1}\aa_{j,1}},A_{j,1})\,\dm\,\cdots\,\dm\,N_2(e^{\sqrt{-1}\aa_{j,r_{j,3}}},A_{j,r_{j,3}})\nn\\
&& \dm\,N_2(e^{\sqrt{-1}\td{\aa}_{j,1}},\td{A}_{j,1})\,\dm\,\cdots\,\dm\,N_2(e^{\sqrt{-1}\td{\aa}_{j,r_{j,4}}},\td{A}_{j,r_{j,4}})\nn\\
&& \dm\,N_2(e^{\sqrt{-1}\bb_{j,1}},B_{j,1})\,\dm\,\cdots\,\dm\,N_2(e^{\sqrt{-1}\bb_{j,r_{j,5}}},B_{j,r_{j,5}})\nn\\
&& \dm\,N_2(e^{\sqrt{-1}\td{\bb}_{j,1}},\td{B}_{j,1})\,\dm\,\cdots\,\dm\,N_2(e^{\sqrt{-1}\td{\bb}_{j,r_{j,6}}},\td{B}_{j,r_{j,6}})\nn\\
&&
\dm\,H(2)^{\dm h_{j,+}}\dm\,H(-2)^{\dm h_{j,-}},\lb{4.9}
\eea
where $\frac{\th_{j,k}}{2\pi}\in\Q\cap(0,1)\bs \{\frac{1}{2}\} $ for $1\le k\le r_{j,1}$, $\frac{\td{\th}_{j,k}}{2\pi}\in(0,1)\bs\Q$ for $1\le k\le r_{j,2}$,
$\frac{\aa_{j,k}}{2\pi}\in\Q\cap(0,1)\bs \{\frac{1}{2}\} $ for $1\le k\le r_{j,3}$, $\frac{\td{\aa}_{j,k}}{2\pi}\in(0,1)\bs\Q$ for $1\le k\le r_{j,4}$,
$\frac{\bb_{j,k}}{2\pi}\in\Q\cap(0,1)\bs \{\frac{1}{2}\} $ for $1\le k\le r_{j,5}$, $\frac{\td{\bb}_{j,k}}{2\pi}\in(0,1)\bs\Q$ for $1\le k\le r_{j,6}$; $N_2(e^{\sqrt{-1}\aa_{j,k}},A_{j,k})$'s and $N_2(e^{\sqrt{-1}\td{\aa}_{j,k}},\td{A}_{j,k})$'s
are nontrivial and $N_2(e^{\sqrt{-1}\bb_{j,k}},B_{j,k})$'s and $N_2(e^{\sqrt{-1}\td{\bb}_{j,k}},\td{B}_{j,k})$'s are trivial, and all non-negative integers satisfy the equality
\bea
p_{j,-}+p_{j,0}+p_{j,+}+q_{j,-}+q_{j,0}+q_{j,+}+r_{j,1}+r_{j,2} + 2\sum_{k=3}^6 r_{j,k}+h_j= n-1. \lb{5.2.0}
\eea

\medskip

{\bf Lemma 4.1.} {\it Under the assumption (F1), for any prime closed geodesic $c_j$, $1\le j\le p$, there holds
	\be
	i(c_j^m)\ge \left[\frac{(2n-3)m}{n-1}\right](n-1),\qquad \forall\ m\in\N,\lb{5.3.0}
	\ee
	and the average index satisfies
	\be
	\hat{i}(c_j)>2n-3. \lb{5.4.0}
	\ee}

{\bf Proof.} By the assumption (F1), since the flag curvature $K$ satisfies $\left(\frac{2n-3}{n-1}\right)^2\left(\frac{\lm}{1 + \lm}\right)^2 < K \le 1$, we can choose $\left(\frac{2n-3}{n-1}\right)^2\left(\frac{\lambda}{\lambda+1}\right)^2<\delta\le K\le 1$. Then by Lemma 2 in \cite{Rad5}, it yields
$$
\hat{i}(c_j)\ge (n-1)\sqrt{\delta}\frac{1+\lambda}{\lambda}>2n-3.
$$

Note that the length of $c_j^m$ satisfies $L(c_j^m)=mL(c_j)\ge m\pi\frac{1+\lambda}{\lambda}>\frac{(2n-3)m}{n-1}\pi/\sqrt{\delta}$ for $m\ge 1$ and $1\le j\le p$ by Theorem 1 of \cite{Rad4}. Then it follows from Lemma 1 of \cite{Rad5} that $i(c_j^m)\ge [\frac{(2n-3)m}{n-1}](n-1)$. \hfill\hb

\medskip

Combining Lemma 4.1 with Theorem 3.5, it follows that
\bea
i(c_j^{m+1})-i(c_j^m)-\nu(c_j^m)\ge i(c_j)-\frac{e(P_{c_j})}{2}\ge 0,\quad\forall\ m\in\N,\ 1\le j\le p.\lb{5.5.0}
\eea
Here the last inequality holds by the fact that $e(P_{c_j})\le 2(n-1)$ and $i(c_j)\ge n-1$.

It follows from (\ref{5.4.0}), Theorem 4.3 in \cite{LoZ} and Theorem 3.6 that for any fixed integer $\bar{m}\in \N$, there exist infinitely many $(p+1)$-tuples $(N, m_1,\cdots, m_p)\in\N^{p+1}$ such that for any $1\le j\le p$ and $1\le m\le \bar{m}$, there holds
\bea i(c_j^{2m_j -m})+\nu(c_j^{2m_j-m})&=&
2N-i(c_j^{m})-\left(2S^+_{P_{c_j}}(1)+2Q_{j}(m)-\nu(c_j^{m})\right), \lb{5.6.0}\\
i(c_j^{2m_j})&\ge& 2N-\frac{e(P_{c_j})}{2},\lb{5.7.0}\\
i(c_j^{2m_j})+\nu(c_j^{2m_j})&\le& 2N+\frac{e(P_{c_j})}{2},\lb{5.8.0}\\
i(c_j^{2m_j+m})&=&2N+i(c_j^{m}),\lb{5.9.0}
\eea
where (\ref{5.7.0}) and (\ref{5.8.0}) follow from (4.32) and (4.33) in Theorem 4.3 in \cite{LoZ} respectively.

Note that by the List 9.1.12 on page 198 of \cite{Lon2}, (\ref{3.21}), (\ref{3.28}), (\ref{4.9}) and $ \nu(c_{j})=p_{j,-}+2p_{j,0}+p_{j,+} $, we have
\bea
&&S^+_{P_{c_j}}(1)=p_{j,-}+p_{j,0},\lb{5.10.0}\\
&&C(P_{c_j})=q_{j,0}+q_{j,+}+r_{j,1}+r_{j,2}+2r_{j,3}+2r_{j,4},\lb{5.11.0}\\
&&Q_{j}(m)=\frac{1+(-1)^m}{2}(q_{j,0}+q_{j,+})+(r_{j,1}+r_{j,3})\nn\\ &&\qquad\qquad\qquad\qquad\qquad\qquad-\sum_{k=1}^{r_{j,1}}\vf\left(\frac{m\theta_{j,k}}{2\pi}\right)-\sum_{k=1}^{r_{j,3}}\vf\left(\frac{m\aa_{j,k}}{2\pi}\right),\lb{5.12.0}\\
&&\nu(c_j^m)=(p_{j,-}+2p_{j,0}+p_{j,+})+\frac{1+(-1)^m}{2}(q_{j,-}+2q_{j,0}+q_{j,+})+2(r_{j,1}+r_{j,3}+r_{j,5})\nn\\
&&\qquad\qquad -2\left(\sum_{k=1}^{r_{j,1}}\vf\left(\frac{m\th_{j,k}}{2\pi}\right)+\sum_{k=1}^{r_{j,3}}\vf\left(\frac{m\aa_{j,k}}{2\pi}\right)
+\sum_{k=1}^{r_{j,5}}\vf\left(\frac{m\bb_{j,k}}{2\pi}\right)\right).\lb{5.13.0}
\eea

By (\ref{5.10.0}), (\ref{5.12.0}) and (\ref{5.13.0}), we obtain
\bea
&&2S^+_{P_{c_j}}(1)+2Q_j(m) -\nu(c_j^{m})\nn\\
&&\qquad=p_{j,-}-p_{j,+}-\frac{1+(-1)^m}{2}(q_{j,-}-q_{j,+})
-2r_{j,5}+2\sum_{k=1}^{r_{j,5}}\vf\left(\frac{m\bb_{j,k}}{2\pi}\right),\nn
\eea
which, together with (\ref{5.6.0}), gives
\bea
i(c_j^{2m_j -m})+\nu(c_j^{2m_j-m})&=&
2N-i(c_j^{m})-p_{j,-}+p_{j,+}+\frac{1+(-1)^m}{2}(q_{j,-}-q_{j,+})\nn\\
&&\quad +2r_{j,5}-2\sum_{k=1}^{r_{j,5}}\vf\left(\frac{m\bb_{j,k}}{2\pi}\right),\quad\forall\ 1\le m\le \bar{m}.\lb{5.14.0}
\eea
By (\ref{5.7.0})-(\ref{5.9.0}), (\ref{5.14.0}), (\ref{5.2.0}), (\ref{5.3.0}) and the fact $e(P_{c_j})\le 2(n-1)$, there holds
\bea
i(c_j^{2m_j-m})+\nu(c_j^{2m_j-m})&\le& 2N+n-1-\left[\frac{(2n-3)m}{n-1}\right](n-1),\nn\\
&&\qquad\qquad\qquad\qquad\forall\ 1\le m\le \bar{m},\lb{5.15.0}\\
2N-(n-1) &\le& 2N-\frac{e(P_{c_j})}{2}\le i(c_j^{2m_j}),\lb{5.16.0}\\
i(c_j^{2m_j})+\nu(c_j^{2m_j})&\le& 2N+\frac{e(P_{c_j})}{2}\le 2N+n-1,\lb{5.17.0}\\
2N+\left[\frac{(2n-3)m}{n-1}\right](n-1)&\le& 2N+i(c_j^{m})=i(c_j^{2m_j+m}), \ \forall\ 1\le m\le \bar{m}.\lb{5.18.0}
\eea

Note that by (\ref{5.5.0}), we have
\bes
&&i(c_j^m)\le i(c_j^{m+1}), \qquad i(c_j^m)+\nu(c_j^m)\le i(c_j^{m+1})+\nu(c_j^{m+1}),\quad\forall m\in \N,
\ees
which, together with (\ref{5.15.0}) and (\ref{5.18.0}), implies
\bea
&&i(c_j^m)+\nu(c_j^m)\le i(c_j^{2m_j -\bar{m}})+\nu(c_j^{2m_j-\bar{m}})\nn\\
&&\qquad\le 2N+n-1-\left[\frac{(2n-3)\bar{m}}{n-1}\right](n-1), \,\forall\ 1\le m\le 2m_j-\bar{m},\lb{5.19.0}\\
&&2N+\left[\frac{(2n-3)\bar{m}}{n-1}\right](n-1)\le i(c_j^{2m_j+\bar{m}})\le i(c_j^{m}),\,\forall\ m\ge 2m_j+\bar{m}.\lb{5.20.0}
\eea

Fix $ \bar{m}\ge 3 $, next we will obtain more accurate index estimates than (\ref{5.15.0}) about $c_j^{2m_j -m}$ for $j\in\{1,\cdots,p\}$ and $m=1,2$.

\medskip

{\bf Lemma 4.2.} {\it Under the assumption (F1), for $j=1,\cdots,p$, we have
	\bea i(c_j^{2m_j -1})+\nu(c_j^{2m_j -1})&\le& 2N-(n-1),\lb{5.23.0}\\
	i(c_j^{2m_j -2})+\nu(c_j^{2m_j -2})&\le& 2N-3(n-1).\lb{5.24.0}\eea}

{\bf Proof.} We give the proofs of (\ref{5.23.0}) and (\ref{5.24.0}) respectively.
\medskip

{\bf Step 1.} {\it The proof of (\ref{5.23.0}).}

By (\ref{3.22}) and (\ref{4.9}), we have
\bea \hat i(c_{j})&=&i(c_{j})+p_{j,-}+p_{j,0}-r_{j,1}-r_{j,2}+\sum_{k=1}^{r_{j,1}}
\frac{\theta_{j,k}}{\pi}+\sum_{k=1}^{r_{j,2}}
\frac{\tilde \theta_{j,k}}{\pi}\nn\\
&\le& i(c_{j})+p_{j,-}+ p_{j,0}+r_{j,1}+r_{j,2}.\lb{5.25.0}
\eea
Combining this with (\ref{5.4.0}), we get
\bea i(c_{j})+p_{j,-}+p_{j,0}+r_{j,1}+r_{j,2}\ge 2n-2. \lb{5.26.0}
\eea
Then, using (\ref{5.14.0}), (\ref{5.26.0}) and (\ref{5.2.0}), we obtain
\bea
i(c_j^{2m_j -1})+\nu(c_j^{2m_j -1})&=&
2N-i(c_j)-p_{j,-}+p_{j,+} \nn\\
&\le & 2N-(2n-2)+p_{j,0}+r_{j,1}+r_{j,2}+p_{j,+}\nn\\
&\le & 2N-(n-1),\lb{5.27.0}
\eea which prove the inequality (\ref{5.23.0}).

\medskip

{\bf Step 2.} {\it The proof of (\ref{5.24.0}).}

\medskip

By (\ref{5.14.0}), we have
\bea
i(c_j^{2m_j-2})+\nu(c_j^{2m_j -2})=
2N-i(c_j^{2})-p_{j,-}+p_{j,+}+q_{j,-}-q_{j,+}.\lb{5.28.0}
\eea
Using (\ref{4.9}), (\ref{3.20}), (\ref{3.22}) and the functions in (\ref{1.3}), we obtain
\bea
i(c_{j}^{2})&=&2(i(c_{j})+p_{j,-}+p_{j,0}-r_{j,1}-r_{j,2})+2\sum_{k=1}^{r_{j,1}} E\left(\frac{\theta_{j,k}}{\pi}\right)+2\sum_{k=1}^{r_{j,2}} E\left(\frac{\tilde{\theta}_{j,k}}{\pi}\right)\nn\\
&&-r_{j,1}-r_{j,2}-p_{j,-}-p_{j,0}-q_{j,0}-q_{j,+},\lb{5.29.0}\\
\hat{i}(c_{j}) &=& i(c_{j})+p_{j,-}+p_{j,0}-r_{j,1}-r_{j,2}+\sum_{k=1}^{r_{j,1}} \frac{\theta_{j,k}}{\pi}+\sum_{k=1}^{r_{j,2}} \frac{\tilde{\theta}_{j,k}}{\pi}.\lb{5.30.0}
\eea
Combining these equalities (\ref{5.29.0}) and (\ref{5.30.0}), we get
\bea
i(c_{j}^{2})&=&2\hat{i}(c_{j})-2\sum_{k=1}^{r_{j,1}} \left\{\frac{\theta_{j,k}}{\pi}\right\}-2\sum_{k=1}^{r_{j,2}} \left\{\frac{\tilde{\theta}_{j,k}}{\pi}\right\}\nn\\
&&+r_{j,1}+r_{j,2}-p_{j,-}-p_{j,0}-q_{j,0}-q_{j,+}.\lb{5.31.0}
\eea
Note that $ \sum_{k=1}^{r_{j,1}} \left\{\frac{\theta_{j,k}}{\pi}\right\}\le r_{j,1}$ and $ \sum_{k=1}^{r_{j,2}} \left\{\frac{\tilde{\theta}_{j,k}}{\pi}\right\}\le r_{j,2}$. Combining this with (\ref{5.28.0}), (\ref{5.31.0}), (\ref{5.4.0}) and (\ref{5.2.0}), we conclude
\bea
i(c_j^{2m_j-2})+\nu(c_j^{2m_j -2})&<&2N-(4n-6)+r_{j,1}+r_{j,2}+p_{j,+}+p_{j,0}+q_{j,0}+q_{j,-}\nn\\
&\le& 2N-(4n-6)+n-1=2N-(3n-5),\lb{5.32.0}
\eea
which implies that
\bea
i(c_j^{2m_j-2})+\nu(c_j^{2m_j -2})\le 2N-(3n-4).\lb{5.33.0}
\eea

Next, assuming that $ i(c_j^{2m_j-2})+\nu(c_j^{2m_j -2})= 2N-(3n-4) $, then by (\ref{5.28.0}) and (\ref{5.32.0}) we have
\bea
i(c_j^{2})+p_{j,-}-p_{j,+}-q_{j,-}+q_{j,+} &=& 3n-4,\lb{5.34.0}\\
r_{j,1}+r_{j,2}+p_{j,+}+p_{j,0}+q_{j,0}+q_{j,-} &=& n-1.\lb{5.35.0}
\eea
Combining (\ref{5.34.0}), (\ref{5.35.0}) and (\ref{5.2.0}), it yields
\bea
p_{j,-}+q_{j,+} &=& 0\lb{5.36.0},\\
i(c_j^{2})-p_{j,+}-q_{j,-} &=& 3n-4\lb{5.37.0}.
\eea
However, from (\ref{5.29.0}), (\ref{5.35.0}) and (\ref{5.36.0}), we deduce
\bea
i(c_j^{2})-p_{j,+}-q_{j,-}&=&r_{j,1}+r_{j,2}+p_{j,-}+p_{j,0}+q_{j,0}+q_{j,+}+p_{j,+}+q_{j,-}\quad(\mod\,2)\nn\\
&=&r_{j,1}+r_{j,2}+p_{j,0}+q_{j,0}+p_{j,+}+q_{j,-}\quad (\mod\,2)\nn\\
&=&n-1\quad (\mod\,2),
\eea
which contradicts (\ref{5.37.0}). Thus $ i(c_j^{2m_j-2})+\nu(c_j^{2m_j -2})\le 2N-(3n-4)-1=2N-3(n-1)$.

\medskip

Now Steps 1-2 complete the proof of Lemma 4.2.\hfill\hb

\medskip

%

The following lemma gives the existence of one elliptic closed geodesic satisfying some special properties on certain positively curved $ (S^{n},F) $.

\medskip

{\bf Lemma 4.3.} (cf. Lemma 3.1 of \cite{Dua1} and Section 3 of \cite{Dua2}) {\it Under the assumption $(\frac{\lm}{1 + \lm})^2 < K \le 1$, there exist at least two prime elliptic closed geodesics on $(S^n,F)$, and one of them, saying $c_{j_{0}}$, satisfying the following properties
	\bea && i(c_{j_{0}}^{2m_{j_{0}}})+\nu(c_{j_{0}}^{2m_{j_{0}}})=2N+n-1,\quad \ol{C}_{2N+n-1}(E,c_{j_{0}}^{2m_{j_{0}}})=\Q,\lb{5.38.0}\\
	&& p_{{j_{0}},-}=q_{{j_{0}},+}=r_{{j_{0}},3}=r_{{j_{0}},4}=r_{{j_{0}},6}=h_{j_{0}}=0,\lb{5.40.0}\\
	&& r_{j_{0},2}=\Delta_1\ge 1, \lb{5.41.0}
	\eea
where the freedom of choose and existence of $N$ and $m_{j_{0}}$ are specified in Theorem 3.6, and all integers in (\ref{5.40.0}) and (\ref{5.41.0}) are defined in Theorem 3.4 and Theorem 3.6.}

\medskip

{\bf Proof of Theorem 1.1.} We use two following steps to prove Theorem 1.1.

\medskip

{\bf Step 1.} {\it The existence of $ n-1 $ prime closed geodesics.}

\medskip

By Lemma 2.6, for every $i\in\N$, there exist some $m(i), j(i)\in\N$ such that
\bea
E( c_{j(i)}^{m(i)})=\kappa_i,\quad
\ol{C}_{2i+\dim(z)-2}(E, c_{j(i)}^{m(i)})\neq 0,\lb{6.1}
\eea
where $\kappa_i\neq\kappa_{i^{\prime}}$ for $ i\neq i^{\prime} $ by Lemma 2.5. It follows from Subsection 2.2 that $\dim(z)=n+1$. By combining Proposition 2.1 with (\ref{6.1}), we have
\bea
i(c_{j(i)}^{m(i)}) \le 2i+n-1 \le i(c_{j(i)}^{m(i)})+\nu(c_{j(i)}^{m(i)}),\qquad \forall\ i\in\N.\lb{6.2}
\eea

Note that for any $i\in G_1\equiv[N-(n-2),N-1]\cap\N$, there holds $2i+n-1\in \bar{G_1}\equiv[2N-(n-3),2N+(n-3)]\cap\N$.

\medskip

{\bf Claim 1.} {\it The following fact holds
\bea
m(i)=2m_{j(i)}, \quad\forall\ i\in G_1. \lb{6.3}
\eea}
\indent In fact, assuming that $ m(i_{0})\neq 2m_{j(i_{0})} $ for some
$ i_{0}\in G_1$. If $ m(i_{0})= 2m_{j(i_{0})}-m $ with $ m\in\N $, then by (\ref{5.15.0}), (\ref{5.19.0}), (\ref{5.23.0}) and $2i_0+n-1\in\bar{G_1}$, we have
$ i(c_{j(i_{0})}^{m(i_{0})})+\nu(c_{j(i_{0})}^{m(i_{0})})\le 2N-(n-1)< 2i_{0}+n-1 $, which contradicts (\ref{6.2}).  If $ m(i_{0})= 2m_{j(i_{0})}+m $ with $ m\in\N $, then by (\ref{5.18.0}), (\ref{5.20.0}) and $2i_0+n-1\in\bar{G_1}$, we have
$ i(c_{j(i_{0})}^{m(i_{0})})\ge 2N+(n-1) > 2i_{0}+n-1 $, which again contradicts (\ref{6.2}). This proves Claim 1.

\medskip

Consequently it yields $ c_{j(i)}\neq c_{j(i^{\prime})}$ for any two different $i,i^{\prime}\in G_1$ since $E(c_{j(i)}^{2m_{j(i)}})=\kappa_i\neq\kappa_{i^{\prime}}=E(c_{j(i^{\prime})}^{2m_{j(i^{\prime})}})$.
Thus there exist at least $n-2$ prime closed geodesics on $(S^{n}, F)$, which is denoted by $\{c_1,\ldots, c_{n-2}\}$ after a permutation of $\{1,\ldots, p\}$, i.e.
\bea
\{j(i)\mid i\in G_1\}=\{1,2,\cdots,n-2\}. \lb{6.3.0}
\eea

By Lemma 4.3, there exists $j_0\in\{1,\dots, p\}$ such that
 $\ol{C}_{2N+n-1}(E, c_{j_0}^{2m_{j_0}})\neq 0$.
So, by Proposition 2.2 and (\ref{5.17.0}), we have
\bea
\ol{C}_{q}(E, c_{j_0}^{2m_{j_0}})=0,\quad \forall q\neq 2N+(n-1). \lb{6.4}
\eea
From (\ref{6.1}) and (\ref{6.4}), we conclude that $c_{j_0}\notin\{c_1,\ldots, c_{n-2}\}$. Hence there exist at least $ n-1 $ closed geodesics on $(S^{n}, F)$.

\medskip

{\bf Step 2.} {\it The existence of another prime closed geodesic.}

\medskip

We carry out this proof by contradiction. Suppose there exist only $n-1$ closed geodesics on $(S^{n}, F)$, i.e. $p=n-1$. And, as inferred in Step 1, we have $c_{j_{0}}=c_{n-1}$.

In order to get a contradiction, we introduce $G_2\equiv [N-(2n-3),N-(n-1)]\cap\N$ and $\bar{G_2}\equiv [2N-3(n-1)+2,2N-(n-1)]\cap\N$. Note that $i\in G_2$ implies that $2i+n-1\in \bar{G_2}$, and $G_1\cap G_2=\bar{G_1}\cap \bar{G_2}=\emptyset$. By (\ref{5.15.0})-(\ref{5.20.0}), (\ref{6.2}), Lemma 4.3 and Proposition 2.2, it follows that $ m(i)=2m_{j(i)}\;\text{or}\;2m_{j(i)}-1$ for any $i\in G_2$.

\medskip

{\bf Claim 2.} {\it We have the following
\bea
m(i)=2m_{j(i)}-1, \quad\forall \; i\in G_2.\lb{6.5}
\eea}
\indent In fact, we assume that $ m(i_{0}^{\prime})=2m_{j(i_{0}^{\prime})} $ for some $i_{0}\in G_2$, then $ j(i_{0}^{\prime})\neq n-1 $ by (\ref{6.1}) and (\ref{6.4}). Thus by Step 1, we must have $ j(i_{0}^{\prime})=j(i_{0}^{\prime\prime}) $ with some $ i_{0}^{\prime\prime}\in\{1,2,\cdots,n-2\} $. Then by (\ref{6.1}) and (\ref{6.3}), we have
$ \kappa_{i_{0}^{\prime}}=E(c_{j(i_{0}^{\prime})}^{2m_{j(i_{0}^{\prime})}})=E(c_{j(i_{0}^{\prime\prime})}^{2m_{j(i_{0}^{\prime\prime})}})=\kappa_{i_{0}^{\prime\prime}} $. This is a contradiction by Lemma 2.5 since $\Nn(S^n,F)<+\infty$. So Claim 2 holds.

\medskip

Consequently it yields $c_{j(i)}\neq c_{j(i^{\prime})}$ for any $i,i^{\prime}\in G_2$ and $i \neq i^{\prime}$ since $E(c_{j(i)}^{2m_{j(i)}-1})=\kappa_i\neq\kappa_{i^{\prime}}=E(c_{j(i^{\prime})}^{2m_{j(i^{\prime})}-1})$. This implies
\bea
\{j(i)\mid i\in G_2\}=\{1,2,\cdots,n-1\}. \lb{6.6}
\eea

\medskip

{\bf Claim 3.} {\it There holds
\bea
\ol{C}_{2N-(n-1)}(E,c_{j}^{m})\neq 0,\;\text{only if}\; j= j(N-(n-1))\; \text{and}\; m=2m_{j(N-n+1)}-1.\lb{6.7}
\eea}
\indent In fact, from (\ref{5.15.0})-(\ref{5.20.0}), Lemma 4.3 and Proposition 2.1, we can deduce that for any $j\in\{1,2,\cdots,n-1\}$ there holds
\bea
\ol{C}_{2N-(n-1)}(E,c_{j}^{m})\neq 0,\quad\text{only if} \ m\in\{2m_{j}, 2m_{j}-1\}.\lb{6.7.1}
\eea

Firstly, assuming that for some $ j_{0}^{\prime}\in\{1,2,\cdots,n-1\} $, there holds
$\ol{C}_{2N-(n-1)}(E,c_{j_{0}^{\prime}}^{2m_{j_{0}^{\prime}}})\\
\neq 0$. This, together with (\ref{5.16.0}), Propositions 2.1 and 2.2, and (\ref{6.4}), implies that
\bea
\ol{C}_{q}(E,c_{j_{0}^{\prime}}^{2m_{j_{0}^{\prime}}})= 0\quad \text{for}\; q\neq 2N-(n-1)\lb{6.7.3}
\eea
and
\bea
c_{j_{0}^{\prime}}\neq c_{j_{0}}, \quad \text{i.e.}\quad j_{0}^{\prime}\neq n-1.\lb{6.7.4}
\eea
By (\ref{6.7.4}) and (\ref{6.3.0}), we have $ j_{0}^{\prime}=j(i_{1}) $ for some $ i_{1}\in G_{1} $, and by (\ref{6.1}) and (\ref{6.3}), we have $ \ol{C}_{2i_{1}+n-1}(E,c_{j(i_{1})}^{2m_{j(i_{1})}})\neq 0 $, i.e., there holds
\bea
\ol{C}_{2i_{1}+n-1}(E,c_{j_{0}^{\prime}}^{2m_{j_{0}^{\prime}}})\neq 0 \quad \text{for some}\;i_{1}\in G_{1}.\lb{6.7.5}
\eea
However, note that $ 2i_{1}+n-1\in \bar{G_1} $, $  2N-(n-1)\in \bar{G_2} $ and $\bar{G_1}\cap\bar{G_2}=\emptyset$, so it yields $2i_{1}+n-1\neq 2N-(n-1)$, which implies that (\ref{6.7.5}) contradicts (\ref{6.7.3}).

Secondly, assuming that for some $ j_{0}^{\prime\prime}\in\{1,2,\cdots,n-1\}\setminus\{j(N-(n-1))\} $, there holds
$
 \ol{C}_{2N-(n-1)}(E,c_{j_{0}^{\prime\prime}}^{2m_{j_{0}^{\prime\prime}}-1})\neq 0,
$
which, together with (\ref{5.23.0}) and Propositions 2.1 and 2.2, implies that
\bea
\ol{C}_{q}(E,c_{j_{0}^{\prime\prime}}^{2m_{j_{0}^{\prime\prime}}-1})= 0\quad \text{for}\; q\neq 2N-(n-1).\lb{6.7.7}
\eea
By (\ref{6.6}), we know that $ j_{0}^{\prime\prime}=j(i_{2}) $ for some $ i_{2}\in G_{2}\setminus\{N-(n-1)\}$, and by (\ref{6.1}) and (\ref{6.5}), we have  $ \ol{C}_{2i_{2}+n-1}(E,c_{j(i_{2})}^{2m_{j(i_{2})}-1})\neq 0 $, then there holds
\bea
\ol{C}_{2i_{2}+n-1}(E,c_{j_{0}^{\prime\prime}}^{2m_{j_{0}^{\prime\prime}}-1})\neq 0 \quad \text{for some}\;i_{2}\in G_{2}\setminus\{N-(n-1)\}.\lb{6.7.8}
\eea
However, note that $ 2i_{2}+n-1\in \bar{G_2}\setminus\{2N-(n-1)\}$, so (\ref{6.7.8}) contradicts (\ref{6.7.7}).

The above two contradictions, together with (\ref{6.7.1}), complete the proof of Claim 3.

\medskip

Combining (\ref{6.7}), (\ref{6.1}), (\ref{6.5}) and Proposition 2.1, we obtain
\bea
M_{2N-(n-1)}=M_{2N-(n-1)}(\ol{\Lm}S^{n},\ol{\Lm}^0 S^{n})=
\dim\ol{C}_{2N-(n-1)}(E, c_{j(N-n+1)}^{2m_{j(N-n+1)}-1})=1.\lb{6.9}
\eea
However, by Theorem 3.6, choosing $ N $ as a multiple of $ n-1 $, we obtain $ 1=M_{2N-(n-1)}\ge b_{2N-(n-1)}=2 $ by Theorems 2.6 and 2.7, which is a contradiction.

\medskip

Now Steps 1 and 2 complete the proof of Theorem 1.1.\hfill\hb

\subsection{Proof of Theorem 1.2 and Theorem 1.3}

In this subsection, we carry out the proofs of Theorems 1.2 and 1.3 under the assumption (F1).

\medskip

{\bf Proof of Theorem 1.2.} Under the assumption (F1), for those closed geodesics $\{c_j\}_{j=1}^p$, if $c_j$ is hyperbolic, then by (\ref{5.16.0}) and (\ref{5.17.0}), there holds
$$i(c_{j}^{2m_j})=2N=i(c_{j}^{2m_j})+\nu(c_{j}^{2m_j}).$$ Thus for every hyperbolic $c_j$, by Proposition 2.1 we have
\bea
\ol{C}_{q}(E, c_{j}^{2m_{j}})=0,\quad \forall\;q\neq 2N.\lb{6.10}
\eea
By (\ref{6.1}), (\ref{6.3}) and (\ref{6.10}), the closed geodesics $\{c_1,\ldots, c_{n-2}\}$ found in Step 1 in the proof of Theorem 1.1 are non-hyperbolic when $ n $ is even, and at least $ n-3 $ of them are non-hyperbolic when $ n $ is odd. By Lemma 4.3, the closed geodesic $c_{j_0}$ is non-hyperbolic. Therefore there are $2[\frac{n}{2}]-1$ non-hyperbolic closed geodesics on $S^{n}$.\hfill\hb

\medskip

{\bf Proof of Theorem 1.3.} Here we mainly follow some ideas in the proof of Theorem 1.4 in \cite{Wan3}. For the reader's conveniences, we give some necessary details.

\medskip

{\bf Claim 4.} {\it There are at least $n-2$ closed geodesics
$c_{j_k}$ for $1\le k\le n-2$ on $(S^{n}, F)$ such that $ c_{j_k}\in \mathcal{V}_\infty(S^{n}, F)$ for $1\le k\le n-2$. }

\medskip

As in the proof of Theorem 1.1, for any $N$ chosen in (\ref{5.6.0})-(\ref{5.9.0})
and $2i+\dim(z)-2\in (2N-n+1,\,2N+n-1)$, there exist some $1\le j(i)\le p$
such that $c_{j(i)}$ is $(2m_{j(i)}, i)$-variationaly visible
by (\ref{6.1}) and (\ref{6.3}). Moreover, if $2N-n+1<i_1\neq i_2<2N+n-1$, then there must hold $j(i_1)\neq j(i_2)$.
Hence the map
\bea
\Psi: (2\N+\dim(z)-2)\cap (2N-n+1, 2N+n-1)&\longrightarrow&\{c_j\}_{1\le j\le p},\nn
\\ 2i+\dim(z)-2&\mapsto& c_{j(i)}\lb{6.11}
\eea
is injective. Here when there are more than one $c_j$ satisfying (\ref{6.1}), we take any one of them. Since the choose of $N$ satisfying (\ref{5.6.0})-(\ref{5.9.0}) is infinite and $\Nn(S^n,F)<+\infty$, so Claim 4 holds.

\medskip

{\bf Claim 5.} {\it Among the $n-2$ closed geodesics $\{c_{j_k}\}_{1\le k\le n-2}$ found in Claim 4, there are at least $n-3$ ones possessing irrational average indices.}

\medskip

By Theorem 3.6, we can choose infinitely many $N$ satisfying (\ref{5.6.0})-(\ref{5.9.0}) and
\be
\frac{N}{\bar{M}\hat i(c_{j})}\in\N\quad{\rm and}\quad\chi_j=0, \qquad{\rm if}\quad \hat i(c_{j})\in\Q.\lb{6.12}
\ee

Now suppose $\hat i(c_{j_{k_1}})\in\Q$ and
$\hat i(c_{j_{k_2}})\in\Q$ hold for some $1\le k_1\neq k_2\le n-2 $.
Then by (\ref{3.29}) and (\ref{6.12}) it yields
\bea
2m_{j_{k_1}}\hat i(c_{j_{k_1}})&=&2\left(\left[\frac{N}{\bar{M}\hat i(c_{j_{k_1}})}\right]+\chi_{j_{k_1}}\right)\bar{M}\hat i(c_{j_{k_1}})\nn\\
&=&2\left(\frac{N}{\bar{M}\hat i(c_{j_{k_1}})}\right)\bar{M}\hat i(c_{j_{k_1}})=2N=2\left(\frac{N}{\bar{M}\hat i(c_{j_{k_2}})}\right)\bar{M}\hat i(c_{j_{k_2}})\nn\\
&=&2\left(\left[\frac{N}{\bar{M}\hat i(c_{j_{k_2}})}\right]+\chi_{j_{k_2}}\right)\bar{M}\hat i(c_{j_{k_2}})
=2m_{j_{k_2}}\hat i(c_{j_{k_2}}).
\lb{6.13}
\eea

In addition, by (\ref{6.11}) we have
\be \Psi(2i_1+\dim(z)-2)=c_{j_{k_1}},\quad\Psi(2i_2+\dim(z)-2)=c_{j_{k_2}}\quad
{\rm for \;some}\; i_1\neq i_2.\lb{6.14}
\ee
Thus by (\ref{6.1}), (\ref{6.3}), Lemma 2.5 and the assumption $\Nn(S^n,F)<+\infty$, we have
\be
E(c_{j_{k_1}}^{2m_{j_{k_1}}})=\kappa_{i_1}\neq\kappa_{i_2}=E(c_{j_{k_2}}^{2m_{j_{k_2}}}).\lb{6.15}
\ee
Since $c_{j_{k_1}}, c_{j_{k_2}}\in\mathcal{V}_\infty(S^{n}, F)$, by Theorem 2.8 there holds
\be
\frac{\hat i(c_{j_{k_1}})}{L(c_{j_{k_1}})}=2\sigma=\frac{\hat i(c_{j_{k_2}})}{L(c_{j_{k_2}})}.\lb{6.16}
\ee
Note that for any closed geodesic $c$ on $(S^{n}, F)$, the following relations hold
\be
L(c^m)=mL(c),\quad \hat i(c^m)=m\hat i(c),\quad L(c)=\sqrt{2E(c)},\qquad\forall m\in\N.\lb{6.17}
\ee

Hence it follows from (\ref{6.15})-(\ref{6.17}) that
\bea
2m_{j_{k_1}}\hat i(c_{j_{k_1}})
&=&2\sigma\cdot 2m_{j_{k_1}}L(c_{j_{k_1}})
=2\sigma\sqrt{2E(c_{j_{k_1}}^{2m_{j_{k_1}}})}
=2\sigma\sqrt{2\kappa_{i_1}}\nn\\
&\neq&2\sigma\sqrt{2\kappa_{i_2}}
=2\sigma\sqrt{2E(c_{j_{k_2}}^{2m_{j_{k_2}}})}
=2\sigma\cdot 2m_{j_{k_2}}L(c_{j_{k_2}})
=2m_{j_{k_2}}\hat i(c_{j_{k_2}}),\nn
\eea
which contradicts (\ref{6.13}), and then there must be
$\hat i(c_{j_{k_1}})\in\R\setminus\Q$ or $\hat i(c_{j_{k_2}})\in\R\setminus\Q$.
Hence there is at most one $1\le k\le n-2$
such that $\hat i(c_{j_{k}})\in\Q$, i.e., there are at least $n-3$ ones possessing irrational average indices.
This proves Claim 5.

\medskip

Suppose $c_{j_k}$ is any closed geodesic found in Claim 5,
then we have $\hat i(c_{j_k})\in\R\setminus\Q$.
Therefore by Theorem 3.4, the linearized Poincar\'e map $P_{c_{j_k}}$
of $c_{j_k}$ must contains a term $R(\theta)$ with $\frac{\theta}{\pi}\notin\Q$.
By Lemma 4.3, the closed geodesic $c_{j_0}$ also has this property.
Moreover, it follows from the proof of Theorem 1.1 that $c_{j_0}\notin \{c_{j_k}\}_{1\le k\le n-2}$.

This completes the proof of Theorem 1.3.\hfill\hb

\subsection{Proof of Theorem 1.4}

To prove Theorem 1.4, in this subsection we assume the following:

\medskip

{\bf (F2)} {\it For any Finsler sphere $(S^{n},F)$ $(n\ge 4)$ satisfying $\frac{9}{4}\left(\frac{\lambda}{\lambda+1}\right)^2<K\le 1$ and $ \lambda<2 $, there holds $\Nn(S^n,F)<+\infty$, denoted by $\{c_j\}_{j=1}^p$.}

\medskip

Firstly, we can estimate the index of closed geodesics on $(S^{n},F)$.

\medskip

{\bf Lemma 4.4.} {\it Under the assumption (F2), for any prime closed geodesic $ c_{j}, 1 \le j \le p $, there holds
	\bea
	i\left(c_{j}^m\right) \geq \left[\frac{3m}{2}\right](n-1) \qquad\forall m \in\N.\lb{4.1}
	\eea}

\medskip

{\bf Proof.} By According to Theorem 1 of \cite{Rad4}, we have
\bea
L(c_{j}^{m})=mL(c_{j})\ge m\pi\frac{\lambda+1}{\lambda}.\lb{4.2}
\eea

By (F2) we can choose $ \delta > \frac{9}{4}\left(\frac{\lambda}{1+\lambda}\right)^2 $
such that $ K \ge \delta$, it together with (\ref{4.2}) gives
$ L(c_{j}^{m})> \frac{3 m\pi}{2\sqrt{\delta}} $, which implies $ i\left(c_{j}^m\right)\geq\left[\frac{3m}{2}\right](n-1) $ by Lemma 1 of \cite{Rad5}.\hfill\hb

\medskip

{\bf Proof of Theorem 1.4.} According to Lemma 4.4, every prime closed geodesic $c_{j}$, $1\le j\le p$, on $(S^{n}, F)$ must satisfy
\bea
i(c_{j}^m)\ge \left[\frac{3m}{2} \right](n-1)\ge n-1\qquad \forall\ m\ge 1, \lb{4.3}
\eea
and
\bea
i(c_{j}^m)\ge \left[\frac{3m}{2} \right](n-1)\ge 3(n-1)\qquad \forall\ m\ge 2. \lb{4.4}
\eea

As in the proof of Theorem 1.1, for every $i\in\N$, there exist some $m(i), j(i)\in\N$ such that
\bea
E( c_{j(i)}^{m(i)})=\kappa_i,\quad
\ol{C}_{2i+n-1}(E, c_{j(i)}^{m(i)})\neq 0,\lb{4.6}
\eea
where $\kappa_i\neq\kappa_{i^{\prime}}$ for $ i\neq i^{\prime} $. And then we have
\bea
i(c_{j(i)}^{m(i)}) \le 2i+n-1 \le i(c_{j(i)}^{m(i)})+\nu(c_{j(i)}^{m(i)})\lb{4.7}
\eea
for any $ i\in\N $. Thus by (\ref{4.3}) and (\ref{4.4}), we obtain
\bea
m(i)=1, \qquad \forall\,\ 1 \le i \le n-2. \lb{4.8}
\eea
Consequently it yields $ c_{j(i)}\neq c_{j(i^{\prime})}$ for $ 1 \le i \neq i^{\prime}\le n-2 $ due to $E(c_{j(i)})=\kappa_i\neq\kappa_{i^{\prime}}=E(c_{j(i^{\prime})})$.
Thus there exist at least $ n-2 $ prime closed geodesics on $(S^{n}, F)$, and we denote these closed geodesics by $\{c_1,\ldots, c_{n-2}\}$ after a permutation of $\{1,\ldots, p\}$.

According to Theorem 2.3, we have $ b_{n-1}=b_{n-1}(\ol{\Lm}S^n,\ol{\Lm}^0 S^n)=1 $. Then, by Theorem 2.4, we have $ M_{n-1}(\ol{\Lm}S^n,\ol{\Lm}^0 S^n)\ge b_{n-1}=1 $, which, together with (\ref{4.4}), implies that there exists $ j_{0}\in\{1,\ldots, p\} $ such that $ \ol{C}_{n-1}(E, c_{j_{0}})\neq 0 $. So, by Proposition 2.2 and (\ref{4.3}), we have
\bea
\ol{C}_{q}(E, c_{j_0})=0,\quad\forall q\neq n-1. \lb{4.9.0}
\eea
From (\ref{4.6}) and (\ref{4.9.0}), we conclude that $c_{j_0}\notin\{c_1,\ldots, c_{n-2}\}$. Hence there exist at least $ n-1 $ closed geodesics on $(S^{n}, F)$.

The proof of Theorem 1.4 is complete.\hfill\hb

\vspace{5mm}






\bibliographystyle{abbrv}

\end{document}